\documentclass[10pt]{article}

\usepackage{amsmath, amssymb, amscd, fancyhdr}

\def\eqref#1{(\ref{#1})}

\def\1{\sqrt{-1}\:}

\newcommand{\cntrct}                
{\hspace{2pt}\raisebox{1pt}{\text{$\lrcorner$}}\hspace{2pt}}

\renewcommand{\bar}{\overline}
\renewcommand{\phi}{\varphi}
\renewcommand{\epsilon}{\varepsilon}

\newcommand{\comment}[1]{{}}

\def\blacksquare{\hbox{\vrule width 4pt height 4pt depth 0pt}}
\def\endproof{\blacksquare}

\makeatletter

\newcounter{Mycounter}[section]
\newcounter{lemma}[section]
\newcounter{claim}[section]
\newcounter{sublemma}[section]
\newcounter{corollary}[section]
\newcounter{theorem}[section]
\newcounter{conjecture}[section]
\newcounter{proposition}[section]
\newcounter{definition}[section]
\newcounter{example}[section]
\newcounter{remark}[section]
\newcounter{problem}[section]
\newcounter{question}[section]
\@addtoreset{equation}{section}
\@addtoreset{footnote}{section}
\makeatother

\begin{document}

\begin{center}
{\LARGE\bf Bounds of some parameters \\[3mm]
of elliptic curve on finite field}
\\[4mm]
Alexey Milovanov
\\[4mm]

{\tt almas239@gmail.com}
\end{center}

{\small
\hspace{0.15\linewidth}
\begin{minipage}[t]{0.7\linewidth}
{\bf Abstract} \\
I prove lower bounds of some parameters of elliptic curve over finite field.
There parameters are closely interrelated with cryptographic stability of elliptic curve.
\end{minipage}
}

\section{Intoduction}

Relations between Diffie-Hellman problem on group points on elliptic curve and Pairing Inversion problem was considered in [1].
It showed the importance of consider the Pairing Inversion problem. No effective algorithm that solve Pairing Inversion problem is known yet.

Algorithm that effective solution the Pairing Inversion problem if some parameter of curve is little was propose in [2], but
heaven knows any curve that has little parameter really.

\section{Diffie-Hellman problem and pairing}

Consider $H_1$, $H_2$ and $H_{r}$ - groups of prime order $r$.
I will write $H_1$ and $H_2$ as additive groups.

Map $e$: $H_1 \times H_2$ $\to$ $H_{r}$ is call
non-degenerate pairing if in fixing any non-identity element
$\in H_i$
$e$ is  isomorphism of group $H_{3-i}$ and $H_r$ ($i \in \{1, 2\}$).

Let we have 2 non-degenerate pairing (they may be same)

$e_1$ and $e_2$
and elements $Y$, $[A]Y$, $[B]Y \in H_1$.
Then we can solve Diffie-Hellman problem (find $[AB]Y$) so:

1)Find $e_1([A]Y, U) = z$ for some non-identity $U$.

2) Inverting pairing $e_1(Y, .) = z$ and find $[A]U$.

3) Find $e_2([B]Y, [A]U) = w$.

4) Inverting $e_2(. ,U) = w$ and find $[AB]Y$.

\section{Notations}

Let $E$ is elliptic curve over $\mathbb{F}_q$.

$r$ is prime number. $G_1: = E(\mathbb{F}_q)[r]$. Let $|G_1| = r$.

Let $k$ is the smallest number that $r | q^k - 1$.

Let $G_r = \{\mu \in \mathbb{F}_{q^k}| \mu^r = 1\} $.

$\pi_q$ - Frobenius endomorphism.

Let $a \in \mathbb{N}$. Write $a$ as $a = \sum_{i=0}^{i = k-1}a_iq^i - \sum_{i=0}^{i = k-1}b_iq^i$ mod$(q^k - 1)$, where $a_i, b_i \in \mathbb{N}_{\ge 0}$

There are many tuples $a_i$, $b_i$ that give valid expression for $a$. We choose so $a_i$ and $b_i$ that
$ S = \sum_{i=0}^{i = k-1}a_i + \sum_{i=0}^{i = k-1}b_i$ is minimal.

{\bf Definition} $D(a) := S$.

\section{About $f^d$-view functions}

{\bf Proposition 1}

For any function $f^d$, where $f \in \bar{\mathbb{F}_q}(E)$ there exists $F\in \bar{\mathbb{F}_q}(E)$
such that: $f^d|_{G_1} = F|_{G_1}$ and
deg$(F) = $deg$(f)D(d)$.

                    {\bf Proof:}
write $d$ as $d = \sum_{i=0}^{i = k-1}a_iq^i - \sum_{i=0}^{i = k-1}b_iq^i$ (mod $q^k - 1$),
where $\sum_{i=0}^{i = k-1}a_i + \sum_{i=0}^{i = k-1}b_i = D(d)$.

 $$ f^d = \prod_{i =0}^{i = k-1} f^{a_iq^i}\prod_{i=0}^{i = k-1} (f^{-1})^{b_iq^i} = \prod_{i =0}^{i = k-1} f^{q^ia_i}\prod_{i=0}^{i = k-1} (f^{-q^i})^{b_i} \eqno(1)$$

Denote by $m_i$ and $l_i$ the function that $f^{q^i} = m_i \circ \pi^i_q$ and
$(f^{-1})^{q^i} = l_i \circ \pi^i_q$.
deg $\pi_q^i = q^i$([3]) , hence deg$(m_i) =$deg$(f) =$ deg$(l_i)$.
Frobenis endomorphism acts identically on points $\in E(\mathbb{F}_q)$ so:

$$ f^{q^i}|_{G_1} = m_i|_{G_1} \text{and} (f^{-1})^{q^i}|_{G_1} = l_i|_{G_1} \eqno(2) $$

Let's denote function:
$F := \prod_{i =0}^{i = k-1} m_i^{a_i}\prod_{i=0}^{i = k-1} l_i^{b_i}$.
$f^d|_{G_1} = F|_{G_1}$ from (1) and (2),
deg$(F) = D(d)$deg$(f)$ from properties $m_i$ and $l_i$.
\endproof

\section{Functions that define pairing}

{\bf Proposition 2}

Let $E$ - is elliptic curve, $G$ - subgroup of $E(\bar{\mathbb{F}_q})$ and $|G| = r$,
$f \in \bar{\mathbb{F}_q}(E)$.
Let
$f^d$ define
non-constant homomorphism
between $G$ and subgroup $\bar{\mathbb{F}_q}$.
Then $d\cdot$deg$(f)\ge (1/6)r$.

{\bf Proof} in [1].
\endproof

In [2] F. Vercauteren proposed algorithm that solve problem inverting pairing in case that deg$(f)D(d)$ is little.

Next statement show that it value can't be little.

\hfill

{\bf Proposition 3}

Let $f \in \bar{\mathbb{F}_q}(E)$ and $f^d$ is isomorphism
between
$G_1$ and $G_r$.
Then $D(d)$deg$(f)\ge (1/6)r$.

{\bf Proof:}
There exist $F \in \bar{\mathbb{F}_q}(E)$:  $f^d|_{G_1} = F|_{G_1}$ и deg$F = D(d)$deg$(f)$ by Proposition 1.
Using Proposition 2 for $f = F$ and $d = 1$ we get required result.
\endproof

\hfill

{\bf Corollary}

Let $f \in \bar{\mathbb{F}_q}(E)$ and $f^d$ is isomorphism
between
$G_1$ and $G_r$.
Let $d = c(q^{k-1}+...q + 1) + d_1$.
Then $D(d_1)$deg$(f)\ge (1/12)r$.

{\bf Proof:}
note that if  $f^{c(q^{k-1}+...q + 1)}f^{d_1}$ is isomorphism between
$G_1$ and $G_r$, then $f^{c(q^{k-1}+...q + 1)(q-1)}f^{d_1(q-1)}$
is isomorphism to
($k > 1$).

$f^{c(q^{k-1}+...q + 1)(q-1)}f^{d_1(q-1)} = f^{d_1(q-1)}$.
Write $d_1$ as $d_1 = \sum_{i=0}^{i = k-1}z_iq^i - \sum_{i=0}^{i = k-1}u_iq^i$ (mod$q^k - 1$),
where $\sum_{i=0}^{i = k-1}z_i + \sum_{i=0}^{i = k-1}u_i = D(d_1)$.

Then $(q-1)d_1 = (z_{k-1} + u_0)q^0 + \sum_{i=1}^{i=k-1}(z_{i-1} + u_i)q^i -
(u_{k-1} + z_i)q^0 - \sum_{i=1}^{i=k-1}(u_{i-1} + z_i)q^i$.
Hence and by definition function $D$ we get that $D((q-1)d_1) \le 2D(d_1)$. Hence and by Proposition 3 we get inequality.
\endproof


\begin{thebibliography}{MM}

\bibitem[1]{S. Galbraith, F. Hess, and F. Vercauteren}
S. Galbraith, F. Hess, and F. Vercauteren
{\em S. Galbraith, F. Hess, and F. Vercauteren
Aspects of pairing inversion, IEEE Transactions on
Information Theory 54 (2008), no. 12, 5719-5728}

\bibitem[2]{F. Vercauteren}
F. Vercauteren
{\em The hidden root problem, Pairing-
Based Cryptography - Pairing, Lecture Notes in
Computer Science, vol. 5209, SpringerLink, Berlin,
2008, pp. 89Ц99}



\bibitem[3]{J. H. Silverman}
J. H. Silverman
{\em The Arithmetic of Elliptic Curves,
Springer-Verlag, GTM 106, 1986. Expanded 2nd Edition, 2009}


\end{thebibliography}
\end{document}